%% file: mcbh-new2.tex
              \documentclass[times,doublespace]{nlaauth}
\usepackage[utf8]{inputenc}
\usepackage[T1]{fontenc}
\usepackage{lmodern}
\usepackage{fixltx2e}
\usepackage{booktabs}
\usepackage{microtype}
\usepackage{tikz}
\usepackage{mathtools}
\usepackage{graphicx}
\usepackage{amsfonts}
\usepackage{color}
\usepackage{multirow}
\usepackage{flafter, float, placeins}
\usepackage{amsmath,amssymb,cancel,units}

\usepackage{pgfplots}
\usepackage{algorithmic} \usepackage[ruled]{algorithm}
\usepackage[unicode=true,plainpages=false]{hyperref}
\hypersetup{colorlinks=true,linkcolor=magenta,anchorcolor=magenta,urlcolor=blue,citecolor=blue}
\newcommand*{\doi}[1]{doi: \href{http://dx.doi.org/#1}{#1}}
\usetikzlibrary{calc,patterns,decorations.pathreplacing,decorations.markings}
\newtheorem{theorem}{Theorem}[section]
\newtheorem{definition}[theorem]{Definition}

\newtheorem{lemma}[theorem]{Lemma}
\def\diam{\mathop{\mathrm{diam}}\nolimits}
\def\dist{\mathop{\mathrm{dist}}\nolimits}

\newcommand{\T}[1]{\mathcal{T}_\mathbf{#1}}
\renewcommand{\t}{\mathbf{t}}
\newcommand{\s}{\mathbf{s}}
\newcommand{\sons}[1]{\bf sons({#1})}
\newcommand{\x}[1]{\chi_\mathbf{#1}}
\renewcommand{\S}[1]{\mathbf{S_{#1}}}
\newcommand{\Rs}[1]{\mathbf{R_{#1}^*}}
\newcommand{\Cs}[1]{\mathbf{C_{#1}^*}}
\newcommand{\pred}[1]{\bf pred({#1})}
\newcommand{\Rp}[1]{\mathbf{R_{#1}^+}}
\newcommand{\Cp}[1]{\mathbf{C_{#1}^+}}
\newcommand{\R}[1]{\mathbf{R_{#1}}}
\renewcommand{\C}[1]{\mathbf{C_{#1}}}
\newcommand{\hx}[1]{\hat{\chi}_\mathbf{#1}}
\newcommand{\p}[1]{\psi_\mathbf{#1}}
\newcommand{\pp}[1]{\psi^p_\mathbf{#1}}
\author{A. Yu. Mikhalev\footnotemark[1] \and I. V. Oseledets\footnotemark[2] \footnotemark[3]}
\date{\today}
\title{}
\hypersetup{
  pdfkeywords={},
  pdfsubject={},
  pdfcreator={Emacs 24.3.50.2 (Org mode 8.2.3c)}}
\begin{document}

\runningheads{A.~Yu.~Mikhalev, I.~V.~Oseledets}{Iterative representing set selection for nested cross approximation}
\title{Iterative representing set selection for nested cross approximation}
\author{A.~Yu.~Mikhalev\affil{1}, I.~V.~Oseledets\affil{1, 2}}
\address{\affilnum{1}
    Skolkovo Institute of Science and Technology,
    Novaya St.~100, Skolkovo, Odintsovsky district, 143025\break
     \affilnum{2}Institute of Numerical Mathematics, Russian Academy of Sciences. Gubkina~St.~8, 119333 Moscow, Russia. \texttt{muxasizhevsk@gmail.com,i.oseledets@skoltech.ru}}
      \cgs{This work was supported by Russian Science Foundation Grant 14-11-00659} 
\begin{abstract}
A new fast algebraic method for obtaining an $\mathcal{H}^2$-approximation
of a matrix from its entries is presented. 
The main idea behind the method is based on the  nested representation
and the maximum-volume principle to select submatrices in low-rank
matrices. A special iterative approach for the computation of
so-called representing sets is established. 
The main advantage of the method is that it uses only the hierarchical
partitioning of the matrix and does not require special ``proxy
surfaces'' to be selected in advance.

The numerical experiments for the electrostatic problem and for the boundary
integral operator confirm the effectiveness and robustness of the
approach. The complexity is linear in the matrix size and polynomial
in the ranks. The algorithm is implemented as an open-source Python
package that is available online.
\end{abstract}
\maketitle
\section{Introduction}
\label{sec-1}

\begin{figure}[H]
\centering
\resizebox{4cm}{!}{
\begin{tikzpicture}[y = -1cm, inner sep = 0, every node/.style = {scale = 1/8}]
\path (0, 0) rectangle +(2, 2) node at +(0, 0) {\usebox{\hmatrixE}};
\end{tikzpicture}}
\caption{Simplest mosaic partitioning. The blocks denoted in white correspond to low-rank matrices.}
\label{fig:fig1}
\end{figure}
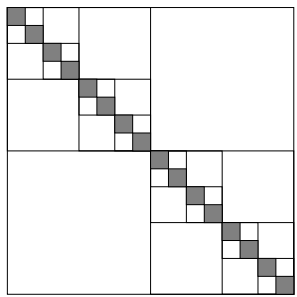

Cross approximation \cite{bebe-cross-2000,tee-cross-2000} is
a basic reduction technique for the
approximation of large low-rank matrices. These matrices
often appear as blocks of dense matrices, coming from the
discretization of non-local operators (in particular, in FEM/BEM applications).
The whole matrix is split into
large blocks related to geometrically separated sets of sources and
receivers, and those blocks are approximated by low-rank
matrices. This partitioning forms the basis for the mosaic-skeleton
method \cite{tee-mosaic-1996} or the
$\mathcal{H}$-matrix format \cite{hack-hmatrix-1999,hk-hmatrix-2000}.
The approximation can be computed only
from the entries of the matrix and from the additional geometrical
information that induces the partititioning of the matrix into
blocks.
The simplest partitioning scheme corresponding to the so-called \emph{weak admissibility condition} \cite{hkk-weak-2004} is presented on
Figure \ref{fig:fig1}, and is usually used for one-dimensional
problems.

\input{include.tex}
In the $\mathcal{H}$-matrix format, the blocks are approximated independently
of each other, which leads to the logarithmic factor in storage and
complexity. The $\mathcal{H}^2$-matrix format \cite{hks-h2matrix-2000} provides a more
effective representation of a matrix based on the analogy with the
well-known fast multipole method (FMM) \cite{gr-fmm-1987}.
The approximation is performed
by taking into account common information about the column and
row bases in the low-rank blocks. This leads to a ``nested''
representation: we do not need to store the factors for each block,
but only the coefficients that allow us to express the bases on the
upper level from the previous level bases (so-called \emph{transfer matrices})
 and coefficients of interactions between bases of sources and receivers (so-called
 \emph{interaction matrices}).
 This structure may lead to significant reduction in memory
and complexity \cite{borm-h2book-2010}.

If the partitioning of a matrix into blocks is fixed, and the
subroutine that allows to compute any prescribed element of a matrix
is given, the main problem is to compute the approximation of such a 
matrix without computing all of its entries (that would give
$\mathcal{O}(N^2)$ complexity since the matrix is dense). For the
$\mathcal{H}$-matrix format,
the problem reduces to the approximation of a low-rank matrix from its
entries. A rank-$r$ matrix can be exactly recovered from $r$ columns
and rows using skeleton decomposition \cite{gtz-psa-1997} and in
the approximate case, the existence of a quasioptimal skeleton
decomposition can be proved based on the maximal volume principle
\cite{gt-maxvol-2001}. In practice, such approximations are computed using
cross approximation techniques \cite{bebe-cross-2000,tee-cross-2000}. 

It is natural to consider a question, whether the
$\mathcal{H}^2$-matrix representation can be computed directly from
the entries of a matrix in a quasi-optimal cost (i.e. the number of
sampled elements is close to the number of parameters in the
representation). This construction can be done efficiently using interpolation \cite{hb-h2interp-2002}, 
but it is not a purely algebraic approach. In \cite{bebe-nested-2012} a \emph{nested cross approximation} was proposed that uses only the entries of the matrix. Two methods were described, ACAGEO and ACAMERGE, which 
rely on different sampling strategies for the \emph{far zone} (in other words, selection of the \emph{representing sets}). This selection was not adaptive for ACAGEO, and for ACAMERGE, as we 
will see in numerical experiments, it may lead to significant loss of accuracy. We propose a new purely algebraic method for the adaptive selection of the representing sets to compute the $\mathcal{H}^2$-matrix approximant  from
the entries. The algorithm was motivated by two completely different
approaches. Its first step comes as a generalization of the classical
Barnes-Hut algorithm \cite{bh-bh-1986}, which can be
considered as one of the first algorithms for the fast approximate computation of dense
matrix-by-vector products. The obtained approximation is satisfactory
for small accuracies (say, $\varepsilon = 10^{-3}$) but does not
perform well when the approximation threshold goes to zero. To make
the method more robust, we establish the analogy between
$\mathcal{H}^2$-matrix approximation and multidimensional tensors, and
apply the ideas from the cross approximation of tensors
\cite{ot-ttcross-2010} to our problem. The top-to-bottom and bottom-to-top tree traversal algorithm, proposed in this paper, is similar to 
the left-right sweeps in the TT-cross algorithm of \cite{ot-ttcross-2010}. 
By numerical experiments we show that the resulting approach is robust in the sense that it is able to achieve 
better accuracy in comparison to the ACAGEO and ACAMERGE algorithms of \cite{bebe-nested-2012}.
 
\emph{Related work.} The problem of computing the $\mathcal{H}^2$-matrix
 approximation from the entries was first considered in this setting in
\cite{bebe-nested-2012}. Our work follows a similar framework but has several
important differences, and the proposed algorithm is more robust
 for high approximation accuracies. Similar approaches for
 constructing the approximation were proposed in the
 \emph{kernel-independent fast multipole methods} \cite{ybz-independent-2004},
but they still use additional information from the problem.
Different wavelet-based methods
\cite{LS-wavebem-1999,amaratunga1994wavelet}, that proved to have
superior performance over hierarchical matrices, suffer from a similar problems: 
they require additional information from the problem and require the computation of Galerkin elements with wavelet basis function.
Some other approaches consider ``hierarchically semiseparable matrices''
\cite{sdc-hss-2007},
this format is very close to $\mathcal{H}^2$.
For this class of matrices (HSS-matrices) randomized algorithm has been 
proposed to recover the representation from the entries
\cite{martinsson-randomhss-2011}.
However, the HSS representation is essentially a one-dimensional $\mathcal{H}^2$-matrix,
and the randomized algorithm of \cite{martinsson-randomhss-2011} 
also requires a fast matrix-by-vector procedure to be available.

We organize our paper as follows. In Section \ref{sec-2} we gather all
necessary mathematical tools to work with $\mathcal{H}^2$-matrices. In
Section \ref{sec-3} we derive the MCBH representation (multicharge
Barnes-Hut) for the $\mathcal{H}^2$-matrix that is main approximation
ansatz, and propose a simple algorithm to compute it.
In Section \ref{sec-4} the connection with good submatrices is
discussed.
In Section \ref{sec-5} the main improvement of the algorithm is proposed. 
In Section \ref{sec-6} we provide numerical experiments.
And, finally, In Section \ref{sec-7} we compare our iterative
approach, proposed in Section \ref{sec-5}, with methods from
\cite{bebe-nested-2012}.

\section{Basic notations and definitions}
\label{sec-2}

A detailed review can be found in the books \cite{hks-h2matrix-2000,borm-h2book-2010}.
Here we give a summary of our notations.

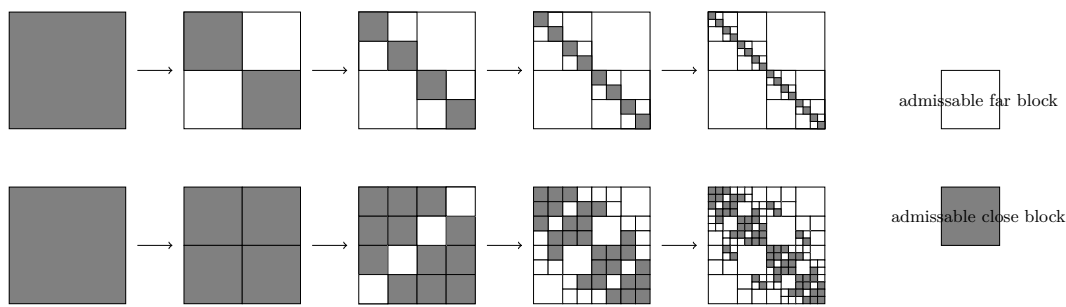
\begin{figure}[H]
\resizebox{14cm}{!}{\input{hmatrix.tex}}
\caption{Examples of mosaic partitioning in one and two dimensions.}
\label{pic:hmatrix}
\end{figure}

To define the $\mathcal{H}^2$-matrix format, we need several standard
definitions. The main idea is to split a matrix into large blocks that
are approximately of low rank.
The matrix element $A_{ij}$ describes interaction between the
$i$-th receiver $x_i$ and $j$-th source $y_j$.
Typical examples are function-related matrices $A_{ij} = f(x_i, y_j)$
and integral operators, discretized by collocation, Nystr\"{o}m or Galerkin
methods with basis functions with local support.
The whole set of receivers will be denoted by indices $\mathbf{I}$ and the
whole set of sources will be denoted by indices $\mathbf{J}$.
For simplicity, we will consider only binary trees, but the extension to
 an arbitrary number of sons is trivial.

\begin{definition}[Row and column cluster trees]
Since sources, denoted by $\mathbf{J}$, are presented by columns and receivers,
 denoted by $\mathbf{I}$, are presented by rows of matrix $A$, let call $\T{I}$
 as row cluster tree and $\T{J}$ as column cluster tree, if

 1) $\mathbf{I}$ is the root of $\T{I}$ and $\mathbf{J}$ is the root of $\T{J}$,

 2) If $\t \in \T{I}$ then $\t$ is a disjoint union of its sons
 $\mathbf{t_1} \in \T{I}$ and $\mathbf{t_2} \in \T{I}$, if $\t \in \T{J}$
 then $\t$ is a disjoint union of its sons $\mathbf{t_1} \in \T{J}$
 and $\mathbf{t_2} \in \T{J}$:
\begin{equation*}
\sons{t} = \{ \mathbf{t_1, t_2} \}.
\end{equation*}
\end{definition}

The cluster trees are constructed for the rows and the columns of the
matrix and are typically based on a certain geometrical construction
(i.e., kd-trees, quadtrees and many others). The construction of the
partition $\mathbf{P}$ of a given matrix into low-rank blocks is based on the
\emph{admissibility condition}. The standard form of the admissibility
condition is also geometrical. A set of receivers $\mathbf{Y}$ and
a set of sources $\mathbf{X}$ are said to be admissible, if
\begin{equation*}
\max\{\diam \mathbf{Y}, \diam \mathbf{X} \} \leq \eta \dist(\mathbf{Y}, \mathbf{X}),
\end{equation*}
where $\eta \ge 0$ (it is equal to 0 for Figures \ref{fig:fig1} and \ref{pic:hmatrix}).
In this case it can be proven that for a class of asymptotically
smooth kernels \cite{tee-cross-2000} the corresponding block can be well-approximated by a
low-rank matrix.
This condition is simple to check, and the recursive partition is
typically generated by a recursive procedure. It starts from the roots
of the cluster trees. If for a given block an admissibility condition
is satisfied, the block $(\t, \s)$ is added to the partition $\mathbf{P}$.
Otherwise, all sons of $\t$ and all sons of $\s$ are processed
recursively. This process creates a list of admissible blocks $(\t, \s)$.

\begin{definition}[$\x{t}, \x{s}$]
Let $\t \in \T{I}$, then $\x{t}$ is a diagonal matrix with the following property:
 $\mathbf{i}$-th diagonal element is 1 if the row $\mathbf{i} \in \t$ and 0 otherwise.
Let $\s \in \T{J}$, then $\x{s}$ is a diagonal matrix with the following property: $\mathbf{i}$-th
diagonal element is 1 if the column $\mathbf{i} \in \s$ and 0 otherwise.
\end{definition}

For any pair of admissible nodes $\t$ and $\s$, their corresponding
submatrix $A(\t, \s)$ can be written as:
\begin{equation*}
A(\t, \s) = \x{t} A \x{s}.
\end{equation*}

\begin{definition}[$\S{t}, \S{s}$]
Assume $\t$ is a node of row cluster tree $\T{I}$.
Then ``set of admissible nodes'', $\S{t} = \{ \mathbf{s_0}, \ldots, \mathbf{s_k} \}$,
 is a set of nodes of a column cluster tree $\T{J}$, so that each pair $(\t, \mathbf{s_i})$
 is admissible if and only if $\mathbf{s_i} \in \S{t}$.
Similarly, $\S{s}$ is a set of admissible nodes for node $\s \in \T{J}$.
\end{definition}

\begin{definition}[$\Rs{t}, \Cs{s}$]
Assume $\t$ is a node of row cluster tree $\T{I}$ with a corresponding set of
 admissible nodes $\S{t}$.
An admissible block row $\Rs{t}$ is defined as:
$$\Rs{t} = \x{t} A \sum_{\s \in \S{t}} \x{s}.$$
Similarly, we define admissible block column $\Cs{s}$ as:
$$\Cs{s} = \left( \sum_{\t \in \S{s}} \x{t} \right) A \x{s}.$$
\end{definition}

\begin{definition}[\pred{t}]
Assume $\t = \mathbf{t_0}$ and $\{ \mathbf{t_1}, \ldots, \mathbf{t_k} \}$ is a set of
 nodes of cluster tree $\T{I}$ or $\T{J}$, such that $\mathbf{t_k}$ is the root of tree,
 and $\mathbf{t_{k-i}}$ is the parent of $\mathbf{t_{k-i-1}}$ for each $\mathbf{i} = 0,
 \ldots, \mathbf{k-1}$.
Then, the set $\{ \mathbf{t_1}, \ldots, \mathbf{t_k} \}$ is called ``predecessors'' of
the node $\t$ and can be written as 
$$\pred{t} = \{ \mathbf{t_1}, \ldots, \mathbf{t_k} \}.$$
\end{definition}

\begin{definition}[$\Rp{t}, \Cp{s}, \R{t}, \C{s}$]
Assume $\t = \mathbf{t_0}$ and $\{ \mathbf{t_1}, \ldots, \mathbf{t_k} \} = \pred{t}$ are predecessors
 of the node $\t \in \T{I}$.
$\{ \S{0}, \ldots, \S{k} \}$ are sets of admissible nodes corresponding to
 $\{ \mathbf{t_0}, \ldots, \mathbf{t_k} \}$.
Then $\Rp{t}$ and $\R{t}$ are defined as:
$$\Rp{t} = \x{t} A \sum_{\mathbf{i} = 1}^{ \mathbf{k} } \sum_{\s \in \S{i}} \x{s},$$
$$\R{t} = \Rp{t} + \Rs{t} = \x{t} A \sum_{ \mathbf{i} = 0}^{ \mathbf{k} } \sum_{\s
 \in \S{i}} \x{s},$$
where $\R{t}$ is called ''block row''.
As for $\mathbf{s_0} = \s \in \T{J}$, we define the following:
$$\{ \mathbf{s_1}, \ldots, \mathbf{s_k} \} = \pred{s},$$
$$\S{i} = \S{s_i},$$
$$\Cp{s} = \left( \sum_{ \mathbf{i} = 1}^{ \mathbf{k} } \sum_{\t \in \S{i}} \x{t} \right) A \x{s} ,$$
$$\C{s} = \Cp{s} + \Cs{s} = \left( \sum_{ \mathbf{i} = 0}^{ \mathbf{k} } \sum_{\t \in \S{i}}
 \x{t} \right) A \x{s},$$
where $\C{s}$ is called ''block column''.
\end{definition}

\begin{figure}[H]
\centering
\input{blockstring.tex}
\caption{Different block rows for node $\t$. $\R{t}, \Rp{t}, \Rs{t}.$}
\label{pic:blockrows}
\end{figure}
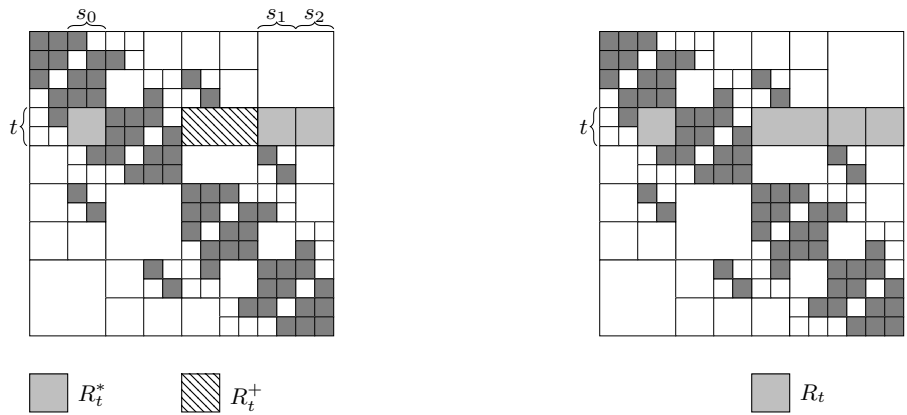

All kinds of block rows are shown on \figurename~\ref{pic:blockrows}.

If $\t \in \T{I}$ is a nonleaf node with $\sons{t} = \{ \mathbf{t_1}, \mathbf{t_2} \}$,
 then block row $\R{t}$ can be computed through block rows of its sons:
\begin{equation*}
\R{t} = \Rp{t_1}+\Rp{t_2}.
\end{equation*}
\section{Nested cross approximation}
\label{sec-3}
In the following, we will always assume exact low-rank property for simplicity, however all 
    the derivation is true for approximate low-rank case as well. In this case, the equality sign should be 
    replaced by approximate equality sign.
The $\mathcal{H}^2$-matrix structure has a simple algebraic 
characterization: all block rows and columns have bounded rank. 
From this property, we can derive a low-rank representation based on the
skeleton decomposition \cite{gtz-maxvol-1997}.
Recall, that the skeleton decomposition for a low-rank matrix $A$ has the following form:
$$A = C \hat{A}^{-1} R,$$
where $C$ consists of \textbf{basis columns} of $A$, $R$ consists of \textbf{basis rows}
of $A$ and $\hat{A}$ is the submatrix of $A$ on the intersection 
of basis rows and columns.
Simple modification of this formula gives us formula for approximating
entire matrix with its basis rows or basis columns:
\begin{equation*}
    A = \tilde{C} R = C \tilde{R},
\end{equation*}
with $\tilde{C} = C \hat{A}^{-1}$ and $\tilde{R} = \hat{A}^{-1} R$.

\begin{figure}[H]
\centering
\begin{tikzpicture}[scale = 0.5]
\draw (0,0) rectangle +(8,8);
\node at (1.5,1.5) {$A$};
\draw[pattern = north east lines] (0,3) rectangle +(8,2);
\draw[pattern = north west lines] (4,0) rectangle +(2,8);
\draw[pattern = north east lines] (10,6) rectangle +(8,2);
\node[right] at (18,7) {$R$};
\draw[pattern = crosshatch] (10,2) rectangle +(2,2);
\node[right] at (12,3) {$\hat{A}$};
\draw[pattern = north west lines] (22,0) rectangle +(2,8);
\node[right] at (24,4) {$C$};
\end{tikzpicture}
\caption{Skeleton decomposition of matrix $A$ by basis rows $R$ and basis columns $C$.}
\label{pic:skeleton}
\end{figure}
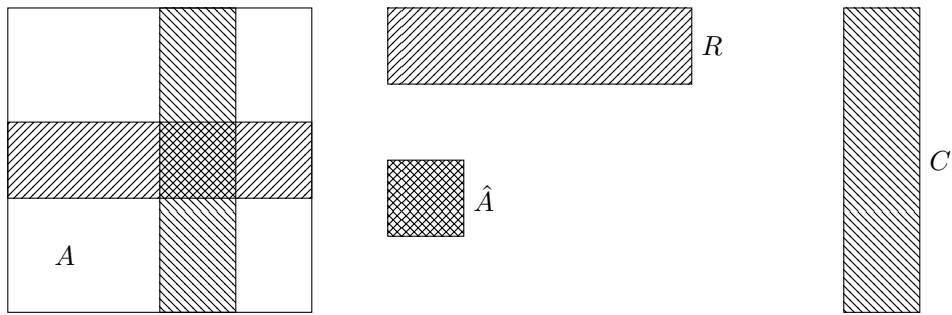

From here, we assume matrix $A$ is a $\mathcal{H}^2$-matrix and every block row and block column is a low-rank matrix and all following equations are written as equalities (for convenience). Since any matrix from real problem can only be close to $\mathcal{H}^2$-matrix (each block column and block row is $\varepsilon$-low-rank), we refer to error analysis from \cite{bebe-nested-2012}.

Let assume, that $\t$ is a nonleaf node of row cluster tree $\T{I}$ with sons $\mathbf{t_1}$
 and $\mathbf{t_2}$.
From the skeleton decomposition for block rows $\R{t}, \R{t_1}$ and $\R{t_2}$ we get the
 following formulas for $\R{t}, \Rp{t_1}$ and $\Rp{t_2}$:
$$\R{t} = U_{\t} \hat{R}_{\t}^{-1} V_{\t},$$
$$\Rp{t_1} = U_{\mathbf{t_1}} (\hat{R}_{\mathbf{t_1}}^{+})^{-1} V_{\mathbf{t_1}},$$
$$\Rp{t_2} = U_{\mathbf{t_2}} (\hat{R}_{\mathbf{t_2}}^{+})^{-1} V_{\mathbf{t_2}},$$
where $V_{\t}, V_\mathbf{t_1}$ and $V_\mathbf{t_2}$ are based on basis rows of
 block rows $\R{t}, \R{t_1}$ and $\R{t_2}$ correspondingly.
Since $\R{t} = \Rp{t_1} + \Rp{t_2}$, block row $\R{t}$ can be
 rewritten as
$$\R{t} = U_{\mathbf{t_1}} \hat{R}_{\mathbf{t_1}}^{-1} V_{\mathbf{t_1}} +
 U_{\mathbf{t_2}} \hat{R}_{\mathbf{t_2}}^{-1} V_{\mathbf{t_2}} = 
 \left[ \begin{matrix} U_{\mathbf{t_1}} (\hat{R}_{\mathbf{t_1}}^{+})^{-1} & U_{\mathbf{t_2}}
 (\hat{R}_{\mathbf{t_2}}^{+})^{-1} \\ \end{matrix} \right]
\left[ \begin{matrix} V_{\mathbf{t_1}} \\ V_{\mathbf{t_2}} \\ \end{matrix} \right].$$
Since $V_\mathbf{t_1}$ and $V_\mathbf{t_2}$ are rows of matrix $\R{t}$, they can be obtained
 via skeleton decomposition of $\R{t}$:
$$ \left[ \begin{matrix} V_{\mathbf{t_1}} \\ V_{\mathbf{t_2}}
\\ \end{matrix} \right] = \hat{U}_{\t} \hat{R}_{\t}^{-1} V_{\t}.$$
Matrix $\hat{U}_{\t} \hat{R}_{\t}^{-1}$ is a
 \emph{local transfer matrix}
 from basis rows of block row $\R{t}$ to basis rows of block
 rows $\R{t_1}$ and $\R{t_2}$.
Here and later we mark it as $M_{\t}$.
Matrices $U_{\mathbf{t_1}} \hat{R}_{\mathbf{t_1}}^{-1}$ and
 $U_{\mathbf{t_2}} \hat{R}_{\mathbf{t_2}}^{-1}$ are \emph{global transfer matrices} from basises of block
 rows $\R{t_1}$ and $\R{t_2}$ to all corresponding rows.
They are marked as $P_{\mathbf{t_1}}$ and $P_{\mathbf{t_2}}$.
Within new notations, we can rewrite low-rank decomposition of block row $\R{t}$:
\begin{equation*}
\R{t} = \begin{bmatrix} P_{\mathbf{t_1}} & P_{\mathbf{t_2}} \end{bmatrix} M_{\t} V_{\t}.
\end{equation*}

For the following argumentation, we will need additional definition:

\begin{definition}[$\hx{t}, \hx{s}$]
Let $\hx{t}$ for $\t \in \T{I}$ be a diagonal matrix with the following property:
$\mathbf{i}$-th diagonal element is 1 if the row $\mathbf{i}$ is in the set of basis
rows for the block row $\R{t}$ and 0 otherwise.
Let $\hx{s}$ for $\s \in \T{J}$ be a diagonal matrix with the following property:
$\mathbf{i}$-th diagonal element is 1 if the column $\mathbf{i}$ is in the set of basis
columns for the block column $\C{s}$ and 0 otherwise.
\end{definition}

For any admissible pair of nodes $(\t, \s)$, corresponding
submatrix $A(\t, \s)$ can be written as follows:
$$A(\t, \s) = \R{t} \x{s} = P_{\t} V_{\t} \x{s} = P_{\t} \hx{t} A
\x{s} = P_{\t} \hx{t} A(\t, \s),$$
where $P_{\t}$ is a global transfer matrix and $V_{\t}$ is a matrix, based on basis rows of $\R{t}$.
Since $\C{s} = U_{\s} P_{\s}$, where $U_{\s}$ is a matrix, based on
basis columns of $\C{s}$, and $P_{\s}$ is a global transfer matrix, we also get
$$A(\t, \s) = A(\t, \s) \hx{s} P_{\s}.$$
Transforming two equations into one:
$$A(\t, \s) = P_{\t} \hx{t} A(\t, \s) \hx{s} P_{\s}.$$
Finally, we get:
\begin{equation*}
A(\t, \s) = P_{\t} \hat{A}(\t, \s) P_{\s},
\end{equation*}
where $\hat{A}(\t, \s)$ is a submatrix, based on basis rows of $\R{t}$
and basis columns of $\C{s}$.
Taking into account that $P_{\t}$ and $P_{\s}$ can be computed recursively, we
get a $\mathcal{H}^2$-type factorization.

First and most simple idea, proposed in the article, is to choose
basis rows for each block row $\R{t}$ hierarchically: use $\hx{t_1}$ and $\hx{t_2}$ to obtain
$\hx{t}$, where $\mathbf{t_1}$ and $\mathbf{t_2}$ are sons of the node
$\t$.
Technically, we can write it as follows:
\begin{equation*}
\hx{t_1} \R{t} + \hx{t_2} \R{t} = M_{\t} \hx{t} \R{t}.
\end{equation*}
Obviously, matrices $M_{\t}$ and $\hx{t}$ can be computed with help of
skeleton decomposition.
The same holds true for basis columns of each block column $\C{s}$:
\begin{equation*}
\C{s} \hx{s_1} + \C{s} \hx{s_2} = \C{s} \hx{s} M_{\s}.
\end{equation*}
Since all basises are chosen hierarchically, this kind of
$\mathcal{H}^2$-decomposition can be named as
 \emph{nested cross approximation}.
In assumption, we generalize proposed proto-method for row cluster tree
$\T{I}$ in Algorithm \ref{alg:prototype}, which can be easily adapted
for column cluster tree $\T{J}$.

\begin{algorithm}[H]
\caption{Computation of basis rows and transfer matrices for each node of row cluster tree}
\label{alg:prototype}
\begin{algorithmic}[1]
\REQUIRE Row cluster tree $\T{I}$, block row $\R{t}$ for each node $\t \in \T{I}$.
\ENSURE nested basises $\hx{t}$ and transfer matrices $M_{\t}$ for each node of $\T{I}$.
\FOR{each $\t \in \T{I}$ \{from bottom of the row cluster tree to top\}}
  \STATE\COMMENT{define auxiliary matrix $\hat{R}$ to compute basis rows and transfer matrix}
  \IF{$\sons{t}$ is empty}
    \STATE $\hat{R} = \R{t}$
  \ELSE[$\sons{t}$ is not empty, $\{ \mathbf{t_1}, \mathbf{t_2} \} = \sons{t}$]
    \STATE{$\hat{R} = \left( \hx{t_1} + \hx{t_2} \right) \R{t}$}
  \ENDIF
  \STATE\COMMENT{compute basis rows $\hx{t}$ and transfer matrix $M_{\t}$}
  \STATE $\hat{R} = M_{\t} \hx{t} \hat{R}$ \COMMENT{obviously, $\hx{t} \hat{R} = \hx{t} \R{t}$}
\ENDFOR
\end{algorithmic}
\end{algorithm}
\section{Maximum volume principle}
\label{sec-4}
Main idea of the nested cross approximation is to select basis rows
and basis columns for special auxiliary matrices ($\hat{R}$ from the
Algorithm \ref{alg:prototype}).
If each block row has precise rank it does not matter what rows we
choose for basises.
Problem arises when block rows are close to low rank matrices.
Skeleton approximation
\begin{equation*}
A \approx C \hat{A}^{-1} R
\end{equation*}
 is about choosing submatrix $\hat{A}$.
If it has very small volume (modulus of determinant), it is very close
 to singular and approximation error is close to infinity.
Guaranteed approximation error can be obtained with help of
 \emph{maximum volume principle} \cite{gtz-psa-1997}: if $\hat{A}$ is
 a submatrix of maximum volume (modulus of determinant) amongst all
 $r \times r$ submatrices, following estimation holds true:
\begin{equation*}
\Vert A - C \hat{A}^{-1} R \Vert_C \le (r+1) \sigma_{r+1},
\end{equation*}
where $\sigma_{r+1}$ is a $r+1$ singular value of matrix $A$ and
$\Vert \cdot \Vert_C$ is an infinite or Chebyshev norm (maximum in
modulus element in matrix).

Finding maximum volume submatrix has exponential complexity, 
 but suboptimal submatrices can be found in a fast way by using greedy algorithms. 
First of all, we can find a good submatrix in a low rank matrix, close to
 $A$:
\begin{equation*}
A \approx U V^T + E,\, A \in \mathbb{C}^{n \times m}, E \in \mathbb{C}^{n \times m}, U \in \mathbb{C}^{n \times r}, V \in \mathbb{C}^{m \times r}.
\end{equation*}
Since 
\begin{equation*}
\det(XY) = \det(X)\det(Y)
\end{equation*}
 for any given square matrices $X$
 and $Y$, computation of good submatrix in $A$
 requires computation of $r$ good rows in matrix $U$ and $r$ good rows
 in matrix $V$.
Recall, that our initial problem was to define basis rows for block
 rows and basis columns for block columns, so we need either left
 factor or right factor.
So without prejudice to the generality, we need to find good
 $r \times r$ submatrix in a $n \times r$ matrix.
Following Algorithm \ref{alg:maxvol}, \emph{maximum-volume algorithm}
 \cite{gostz-maxvol-2010} (which we call \emph{maxvol}), solves this
 problem.

Suppose we need to find a good $r \times r$
  submatrix in a matrix
  $A \in\mathbb{C}^{n \times r}$, with $n>r$.
 First step of \emph{maxvol} algorithm is to find a nonsingular
 submatrix of matrix $A$, so we use the LU decomposition with row
 pivoting, which requires $O(nr^2)$ operations. In other words, initialization
  gives us a submatrix $\hat{A} \in \mathbb{C}^{r \times r}$ and
  coefficients $C \in \mathbb{C}^{n\times r}$ such that $A=C\hat{A}$. Each iteration of the
  algorithm swaps two rows in order to increase the volume of the submatrix. 
  As shown in \cite{gostz-maxvol-2010} each
 iteration requires only $O(nr)$ operations since it finds absolute
 maximum element in $C$ and applies rank-1 update to it.
 In practice the algorithm converges very fast, thus 
 So, we can estimate the complexity of \emph{maxvol} algorithm as $O(nr^2)$ operations.

\begin{algorithm}[H]
\caption{\emph{maxvol} \cite{gostz-maxvol-2010} algorithm}
\label{alg:maxvol}
\begin{algorithmic}[1]
\REQUIRE Nonsingular matrix $A \in \mathbb{C}^{n \times r}, n > r$
\ENSURE Good submatrix $\hat{A} \in \mathbb{C}^{r \times r}$ and
coefficients $C \in \mathbb{C}^{n \times r}$ such, that $A = C \hat{A}$
\STATE Find nonsingular $r \times r$ submatrix $\hat{A}$ in matrix $A$
\STATE\COMMENT{i.e. by pivots from LU factorization}
\STATE $C \leftarrow A \hat{A}^{-1}, \{i,j\} \leftarrow
\mathrm{argmax}(\vert C\vert)$
\STATE\COMMENT{$\vert C\vert$ means
  element-wise modulus, $i$ and $j$ are row and column numbers of
  maximum in modulus element in $C$}
\WHILE{$\vert C_{ij}\vert > 1$}
\STATE\COMMENT{$C_{ij}$ is element on intersection of $i$-th row and
  $j$-th column of $C$}
\STATE $\hat{A}_j \leftarrow A_i$
\STATE\COMMENT{$A_i$ is $i$-th row of $A$, $\hat{A}_j$ is $j$-th row of $\hat{A}$}
\STATE $C \leftarrow C-C_j(C_i-e_j)/C_{ij}$
\STATE\COMMENT{$C_i$ is $i$-th
  row and $C_j$ is $j$-th column of $C$, $e_j$ is $j$-th row of
  identity matrix of size $r \times r$}
\STATE $\{i,j\} \leftarrow \mathrm{argmax}(\left|C\right|)$
\ENDWHILE
\RETURN $C, \hat{A}$
\end{algorithmic}
\end{algorithm}
\section{Practical algorithm for nested cross approximation}
\label{sec-5}
Main disadvantage of the Algorithm \ref{alg:prototype} is
computational cost.
Auxiliary matrix $\hat{R}$ of each leaf node $\t \in \T{I}$ contains
$n_{\t}$ (number of receivers, corresponding to node $\t$) nonzero rows and
$N_{\t}$ (number of sources in far zone of $\t$ or $\pred{t}$) nonzero columns.
Since $\t$ is a nonleaf node, $N_{\t}$ is close to number of sources.
Obviously, minimum number of operations to compute basis rows for node
 $\t$ is a multiplication of $n_{\t}$ and $N_{\t}$.
Summing it for all leaf nodes we get a $O(NM)$ complexity, where $N$
 is number of sources and $M$ is number of receivers.
So we need to reduce the number of nonzero columns of $\hat{R}$, i.e. to
 select a certain \emph{representing set}.

Let assume that we already have representing sets for each node of
$\T{I}$ and $\T{J}$.
It is again convenient to use diagonal matrices to work with
representing sets.
\begin{definition}[$\p{t}, \p{s}$]
Let $\p{t}$ be a diagonal matrix with the following property:
$\mathbf{i}$-th diagonal element is 1 if the column $\mathbf{i}$ is in the representing set for the block row $\R{t}$ and 0 otherwise.
Let $\p{s}$ be a diagonal matrix with the following property:
$\mathbf{i}$-th diagonal element is 1 if the row $\mathbf{i}$ is in the representing set for the block column $\C{s}$ and 0 otherwise.
\end{definition}

If a good representing set for each block row $\R{t}$ is already known,
 then basis rows and transfer matrix for each block row can be
 computed with help of a small matrix $\hat{R} \p{t}$ instead of
 $\hat{R}$.
\emph{Maxvol} algorithm \ref{alg:maxvol} can be used to find basis
 rows or columns together with transfer matrices in an efficient way.
Finally, we get a prototype Algorithm \ref{alg:precomputed}.

\begin{algorithm}[H]
\caption{Computation of basis rows and transfer matrices for row cluster tree with given representing sets.}
\label{alg:precomputed}
\begin{algorithmic}
\REQUIRE Row cluster tree $\T{I}$, block row $\R{t}$ and representing set $\p{t}$ for each node $\t \in \T{I}$, accuracy parameter $\varepsilon$.
\ENSURE basis rows $\hx{t}$ and transfer matrix $M_{\t}$ for each block row $\R{t}$.
\FOR{$\t \in \T{I}$ \{from bottom of the row cluster tree to top\}}
  \STATE\COMMENT{define auxiliary matrix $\hat{R}$ to compute basis rows and transfer matrix}
  \IF{$\sons{t}$ is empty}
    \STATE{$\hat{R} = \R{t} \p{t}$}
  \ELSE[$\sons{t}$ is not empty, $\{ \mathbf{t_1}, \mathbf{t_2} \} = \sons{t}$]
    \STATE{$\hat{R} = \left(\hx{t_1} + \hx{t_2} \right) \R{t} \p{t}$}
  \ENDIF
  \STATE\COMMENT{truncated svd + maxvol}
  \STATE{$U, S, V = svd(\hat{R}, tol = \varepsilon); M_{\t}, \hx{t} = maxvol(U)$}
\ENDFOR
\end{algorithmic}
\end{algorithm}

\subsection{Partially fixed representing sets}
\label{sec-5-1}

The main problem is how to select these representing sets adaptively.
Previous approaches \cite{bebe-nested-2012,ybz-independent-2004,hg-directsolver-2012} used geometrical constructions.
In this paper we propose a purely algebraical method: each representing
 set is divided into self and predecessors parts, each of which is
 calculated separately.
First of all, assume that predecessors part of each representing set is
 predefined, then to find rather good self part of representing set for
 a node $\t$, we can use basis rows or columns of all nodes in $\S{t}$
 or $\sons{\S{t}}$ if these basises are already defined.
This leads us to a \emph{level-by-level} algorithm: start on the
 bottom of trees $\T{I}$ and $\T{J}$, calculate basis columns for each
 node of the $\T{J}$ on current level of the tree, calculate basis
 rows for each node of the $\T{I}$ on the same level and then go up
 one level.
Under assumption of fixed predecessors part $\pp{t}$ of each representing set,
 we summarize \emph{level-by-level} idea in Algorithm \ref{alg:mcbh-upward}.

\begin{algorithm}
\caption{Computation of basis rows/columns and transfer matrices for each node of rows/columns cluster tree with given predecessors part of each representing set.}
\label{alg:mcbh-upward}
\begin{algorithmic}[1]
\REQUIRE Cluster trees $\T{I}$ and $\T{J}$, block row/column $\R{t}$ ($\C{t}$) and predecessors part $\pp{t}$ of representing set for each node $\t$, accuracy $\varepsilon$
\ENSURE basis rows/columns $\hx{t}$ and transfer matrix $M_{\t}$ for each node $t \in \T{I}$ or $\t \in \T{J}$
\STATE\COMMENT{Assume, that both $\T{I}$ and $\T{J}$ have $\mathbf{level\_count}$ levels}
\FOR{$current\_level = \mathbf{level\_count}$ \TO 1 \{from bottom of trees to top\}}
 \FOR{$\s \in \T{J}$ \{on the $current\_level$\}}
  \STATE{$\hat{\psi} = \pp{s}$ \{initialize $\hat{\psi}$ with predecessors part\}}
  \FOR{$\t \in \S{s}$}
   \IF{$\sons{t}$ is empty}
    \STATE{$\hat{\psi} += \hx{t}$ or $\hat{\psi} += \x{t}$ (whether $\hx{t}$ is defined or not)}
   \ELSE[$\sons{t}$ is not empty, $\{ \mathbf{t_1}, \mathbf{t_2} \} = \sons{t}$]
    \STATE{$\hat{\psi} += \hx{t_1} + \hx{t_2}$}
   \ENDIF
  \ENDFOR
  \IF{$\sons{s}$ is empty}
   \STATE{$\hat{R} = \hat{\psi} \C{s}$, in other terms $\hat{R} = \hat{\psi} A \x{s}$}
  \ELSE[$\sons{s}$ is not empty, $\{ \mathbf{s_1}, \mathbf{s_2} \} = \sons{s}$]
   \STATE{$\hat{R} = \hat{\psi} \C{s} \left( \hx{s_1}+\hx{s_2} \right)$, in other terms $\hat{R} = \hat{\psi} A \left( \hx{s_1}+\hx{s_2} \right)$}
  \ENDIF
  \STATE{$U, S, V = svd(\hat{R}, tol = \varepsilon); M_{\s}, \hx{s} = maxvol(V)$}\COMMENT{truncated svd + maxvol}
 \ENDFOR
 \FOR{$\t \in \T{I}$ \{on the $current\_level$\}}
  \STATE{$\hat{\psi} = \pp{t}$ \{initialize $\hat{\psi}$ with predecessors part\}}
  \STATE{$\hat{\psi} += \sum_{\s \in \S{t}} \hx{s}$}
  \IF{$\sons{t}$ is empty}
   \STATE{$\hat{R} = \R{t} \hat{\psi}$, in other terms $\hat{R} = \x{t} A \hat{\psi}$}
  \ELSE[$\sons{t}$ is not empty, $\{ \mathbf{t_1}, \mathbf{t_2} \} = \sons{t}$]
   \STATE{$\hat{R} = \left( \hx{t_1} + \hx{t_2} \right) \R{t} \hat{\psi}$, in other terms $\hat{R} = \left( \hx{t_1} + \hx{t_2} \right) A \hat{\psi}$}
  \ENDIF
  \STATE{$U, S, V = svd(\hat{R}, tol = \varepsilon); M_{\t}, \hx{t} = maxvol(U)$}\COMMENT{truncated svd + maxvol}
 \ENDFOR
\ENDFOR
\end{algorithmic}
\end{algorithm}
\subsection{Representing sets and iterations}
\label{sec-5-2}

How to get predecessors part $\pp{t}$ of each representing set?
Let assume that all basis rows and columns are given.
Obviously, union of basis columns of nodes from
 $\S{t} \cup \S{\pred{t}}$ is a good representing set for block
 row $\R{t}$.
If nonleaf node $\t$ has sons
$\{ \mathbf{t_1}, \mathbf{t_2} \} = \sons{t}$,
 then predecessors part of representing set for node $\mathbf{t_1}$
 can be set equal to just full representing set of $\t$.
This simple idea works perfectly with our previous approach:
\begin{enumerate}
\item with a given predecessors part of representing sets, compute basis
rows and columns,
\item with given basises, update representing sets,
\item shift newly calculated representing sets to produce predecessors part of
representing sets
\item go to step 1.
\end{enumerate}
Obviously, we get an alternating iterative method of computing basises
 and representings.
If we recalculate representings just as union of corresponding basises,
 size of each representing set will depend on number of nodes in
 $\S{t} \cup \S{\pred{t}}$.
This problem can be easily solved by additional resampling in
 top-to-bottom fashion, proposed in Algorithm \ref{alg:mcbh-downward}.

\begin{algorithm}
\caption{Computation of full representing sets with given basises}
\label{alg:mcbh-downward}
\begin{algorithmic}[1]
\REQUIRE Cluster trees $\T{I}$ and $\T{J}$, basis rows/columns $\hx{t}$ for each node $\t$
\ENSURE Representing set $\p{t}$ for each node $\t$
\STATE\COMMENT{Assume, that both $\T{I}$ and $\T{J}$ have $\mathbf{level\_count}$ levels}
\FOR{$current\_level = 1$ \TO $\mathbf{level\_count}$\{from bottom of trees to top\}}
 \FOR{$\t \in \T{I}$ \{on the $current\_level$\}}
  \STATE\COMMENT{Initialize $\hat{\psi}$ as a resampled union of basis rows/columns from $\S{\pred{t}}$}
  \IF{$current\_level$ is 1 \{check if node $\t$ is the root node\}}
   \STATE{$\hat{\psi} = 0$}
  \ELSE[$\mathbf{p}$ is the parent node of $\t$]
   \STATE{$\hat{\psi} = \p{p}$}
  \ENDIF
  \STATE{$\hat{\psi} += \sum_{\s \in \S{t}} \hx{s}$}\COMMENT{Add basis rows/columns from $\S{t}$}
  \STATE{$\p{t} = maxvol(\hx{t} A \hat{\psi})$}\COMMENT{resample $\hat{\psi}$ into $\p{t}$}
 \ENDFOR
\ENDFOR
\end{algorithmic}
\end{algorithm}

\subsection{Iterative Multicharge method}
\label{sec-5-3}

Finally, the proposed MCBH algorithm is an iterative method and is summarized in
Algorithm \ref{alg:total}.
Initialization is very easy: basis rows and columns are empty.
Each iteration consists of 3 steps:
\begin{enumerate}
\item Obtaining representing sets with Algorithm \ref{alg:mcbh-downward}
for each cluster tree,
\item Shifting representing sets into predecessors part of representing sets,
\item Computing basis rows/columns and transfer matrices with Algorithm
   \ref{alg:mcbh-upward}.
\end{enumerate}

Each iteration uses all the information of the previous iteration and increases accuracy.
Numerical examples showed, that 1-2 iterations are usually enough to
get the desired accuracy.

\begin{algorithm}
\caption{Iterative MCBH algorithm.}
\label{alg:total}
\begin{algorithmic}
\STATE Initialize basis rows/columns as empty
\WHILE{did not get required accuracy}
\STATE run Algorithm \ref{alg:mcbh-downward} for each cluster tree to get representing sets
\STATE reevaluate representing sets into predecessors part of representing sets
\STATE run Algorithm \ref{alg:mcbh-upward} to obtain new improved basis rows/columns and transfer matrices
\ENDWHILE
\end{algorithmic}
\end{algorithm}
\subsection{Complexity estimates}
\label{sec-5-4}
    Since the form of $\mathcal{H}^2$-factorization presented in \ref{sec-3} 
    is the same, as
in \cite{bebe-nested-2012}, storage and error analysis holds the same
as in \cite{bebe-nested-2012}.
Only thing to be analised is complexity.
Main approach is formulated as step-by-step
algorithm \ref{alg:total}, so we analyse each step separately.

Let assume we have a matrix $A \in \mathbb{C}^{N \times M}$, that is already partitioned into blocks with row cluster tree $\T{I}$ and column cluster tree $\T{J}$. For simplicity, assume that all the ranks are equal to $r$. 
Then we have to do the following steps.\\
Step 1: Initialize basis rows/columns as empty.
Obviously, it requires $O(\left|\T{I}\right|+\left|\T{J}\right|)$ operations ($\left|\mathcal{T}\right|$ is a number of nodes in cluster tree $\mathcal{T}$).\\
Step 2: run Algorithm \ref{alg:mcbh-downward} for each cluster tree to get representing sets.
Algorithm \ref{alg:mcbh-downward} works in a top-to-bottom manner with each node of clusters $\T{I}$ and $\T{J}$.
On the zero iteration (we enumerate iterations from 0), basises are empty, so this step does not require any operations.
Assume, that for any other iteration basis size for each node of cluster trees $\T{I}$ and $\T{J}$ is equal to constant $r$.
Also, assume any node $t \in \T{I} \cup \T{J}$ has $C_F$ admissibly far nodes and any leaf node $t \in \T{I} \cup \T{J}$ has $C_C$ admissibly close nodes.
Then, computation of representing sets with given basises for all nodes will require $O(maxvol((C_F+1)r \times r))*(\left|\T{I}\right|+\left|\T{J}\right|) = O(C_Fr^3(\left|\T{I}\right|+\left|\T{J}\right|))$ operations.\\
Step 3: reevaluate representing sets into the predecessors part of representing sets. This step is a simple move of representing sets from parent to child for each node (except leaves) of cluster trees $\T{I}$ and $\T{J}$.
Obviously, it requires $O(\left|\T{I}\right|+\left|\T{J}\right|)$ operations.\\
Step 4: run Algorithm \ref{alg:mcbh-upward} to obtain new improved basis rows/columns and transfer matrices.
Assume each non-leaf node has $K$ children, each leaf node corresponds to $Kr$ rows or columns.
Then, SVD reduction of Algorithm \ref{alg:mcbh-upward} works with $Kr \times (C_FKr+r)$ matrices for each node $t \in \T{J}$ and with $Kr \times (C_Fr+r)$ for each node $t \in \T{I}$.
Since $K \ge 2$, we can estimate the complexity of all SVD reductions in Algorithm \ref{alg:mcbh-upward} as $O(\left|\T{J}\right| C_FK^3r^3)$ operations.
Each \emph{maxvol} procedure works with a $Kr \times r$ matrix, so, totally for all nodes, it requires $O((\left|\T{I}\right|+\left|\T{J}\right|)Kr^3)$ operations.\\
With further simplifications ($N = M$, complexity of SVD is of the same order as the complexity of \emph{maxvol}), total complexity is about $O(N_{iters}C_FK^3r^3 \left|\T{J}\right|)$ operations.
As we already assumed each leaf node has exactly $Kr$ elements in it, so we have $\left|\T{J}\right| = O(MK^{-1}r^{-1}) = O(NK^{-1}r^{-1})$.
So, the complexity of the approximation of far interactions in the matrix $A$ requires $O(C_FNK^2r^2)$ operations.
Since number of leaf nodes is $O(NK^{-1}r^{-1})$ and each leaf node requires $O(C_CK^2r^2)$ to save close interactions, we get following complexity on close interactions: $O(C_CNKr)$ operations.
Totally, complexity of proposed MCBH algorithm \ref{alg:total} is $O(C_FNK^2r^2+C_CNKr)$ operations.
where $C_F$ is a mean number of admissibly far nodes for all nodes, $C_C$ is a mean number of admissibly close nodes for all leaf nodes, $N$ is a number of rows and columns, $K$ is a number of children nodes for each node of cluster trees and $r$ is a mean basis size (rank) of each node.

Another interesting question is the number of matrix entries, required to compute approximation with Algorithm \ref{alg:total}.
Often, each matrix element is computed via complex function, which require many operations just for a single matrix element.
Matrix elements are only used in close interactions and steps 2 and 4 of Algorithm \ref{alg:total}.\\
Close interactions: storage cost is the same, as number of operations ($O(C_CNKr)$).
Step 2: run Algorithm \ref{alg:mcbh-downward} for each cluster tree to get representing sets.
This operation requires computation of \emph{maxvol} for $(C_F+1)r \times r$ submatrices of $A$ for each node.\\
Step 4: run Algorithm \ref{alg:mcbh-upward} to obtain new improved basis rows/columns and transfer matrices.
Each iteration uses $Kr \times (C_FKr+r)$ submatrix of $A$ for each node $t \in \T{J}$ and $Kr \times (C_Fr+r)$ submatrix of $A$ for each node $t \in T{I}$.\\
So, each iteration uses $O(C_FK^2r^2\left|\T{J}\right|)$ matrix entries.
Assuming $\left|\T{J}\right| = O(NK^{-1}r^{-1})$, we summarize total number of matrix elements used: $N_{iter}$ iterations of MCBH require $O(N_{iter}C_FNKr+C_CNKr)$ matrix values.

\section{Numerical experiments}
\label{sec-6}
This section contains two numerical examples.
The first is the following:
compute electric field potential, created by $n$ given particles $X_1, \ldots, X_n$ with charges $q_1, \ldots, q_n$, in coordinates of particles $X_1, \ldots, X_n$.
Particles are distributed randomly (uniform distribution) over a cube
$[0, 1]^3$.
This problem can be reformulated as computation of a matrix-vector product:
\begin{equation}
\label{eq:nbody}
Aq=f,
\end{equation}
where $q$ is a vector of charges of particles, $f$ is desired electric
field potential and matrix $A$ is defined as follows:
\begin{equation}
\label{eq:nbodyaij}
A_{ij} = \left\{ \begin{array}{ll}\frac{1}{\left| X_i-X_j \right|}, & i \ne j \\ 0, & i = j \end{array} \right.
\end{equation}

The second example is a solvation problem in the framework of polarized continuum model
 \cite{tp-continuus-1994,ct-implicit-1999,ta-rapidpcm-2001}.
It arises in computer modeling of drugs:
 find the surface charge density $\sigma$ on a given solvent excluded
 surface $\Omega$, such that
\begin{equation}
\label{eq:sigmar}
\sigma(\mathbf{r}) = \frac{1-\varepsilon}{2 \pi
(1+\varepsilon)}\left( \sum_i \frac{Q_i ((\mathbf{r-R_i})\cdot \mathbf{n})}{\left|
\mathbf{r-R_i} \right|^3} + \int_{\Omega}
\frac{\sigma(\mathbf{r'}) ((\mathbf{r-r'}) \cdot \mathbf{n})}{\left|
\mathbf{r-r'} \right|^3} dS' \right),
\end{equation}
where $Q_i$ is a charge of $i$-th atom in molecule,
$\mathbf{R_i}$ is a position vector of $i$-th atom in molecule,
$\mathbf{r}$ is a radius vector to surface,
$\mathbf{n}$ is a perpendicular from surface to solvent
and $\varepsilon$ is a relative permittivity.
The surface is approximated by discrete elements with the Nystr\"{o}m method for
the off-diagonal elements and the diagonal elements are computed from
the identities:
\begin{equation}\label{intomega}
\begin{split}
& \int_{\Omega}\frac{((\mathbf{r-r'}) \cdot \mathbf{n'})}{\left|
\mathbf{r-r'} \right|^3} dS' = 2\pi, \\
& \int_{\Omega_{\varepsilon}} \frac{\sigma(\mathbf{r'}) ((\mathbf{r-r'}) \cdot \mathbf{n})}{\left|
\mathbf{r-r'} \right|^3} dS' \approx \sigma(\mathbf{r}) \left( 2 \pi - \int_{\Omega
\setminus \Omega_{\varepsilon}} \frac{((\mathbf{r-r'}) \cdot
\mathbf{n'})}{\left| \mathbf{r-r'} \right|^3} dS' \right).
\end{split}
\end{equation}
So, after discretizing and integrating over each discrete element, we get following system:
\begin{equation*}
Aq=BQ,
\end{equation*}
where $q$ is a vector of charges of discrete elements and $Q$ is a vector of charges of atoms in molecule.
Matrices $A$ and $B$ are defined as follows:
\begin{equation}
\label{eq:aij}
A_{ij} = \left\{ \begin{array}{ll}
\frac{(\varepsilon-1)}{4\pi(1+\varepsilon)} \frac{((\mathbf{r_i-r_j})
\cdot \mathbf{n_i}) S_i}{\left| \mathbf{r_i-r_j} \right|^3}, & i \ne j
\\ \frac{\varepsilon}{1+\varepsilon} -\sum_{k \ne j} A_{kj}, & i = j
\end{array} \right.,
\end{equation}
\begin{equation*}
B_{ij} = \frac{1-\varepsilon}{4\pi(1+\varepsilon)}
\frac{((\mathbf{r_i-R_j}) \cdot \mathbf{n_i}) S_i}{\left|
\mathbf{r_i-R_j} \right|^3},
\end{equation*}
where $\mathbf{r_i}$ is a radius vector to center of $i$-th discrete element,
$\mathbf{n_i}$ is a perpendecular to $i$-th discrete element
and $S_i$ is an area of $i$-th discrete element.

In the following numerical examples we study only the approximation of
the matrix \eqref{eq:aij}; we are not considering the problem of solution of linear systems with such matrices.
\subsection{Implementation remarks}
\label{sec-6-1}
MCBH algorithm \ref{alg:total} is implemented in Python with the most
computationally intensive parts reimplemented in Cython
[\url{http://cython.org}]. 
For the basic linear algebra tasks the MKL
library is used.
The code can be obtained at \\ \url{http://bitbucket.org/muxas/h2tools}.
To build the hierarchical tree, we used recursive
inertial bisection with admissibility parameter $\eta = 0$ and \emph{block\_size} = 50, where \emph{block\_size} is a constant, such that we divide cluster into subclusters only if it has more than \emph{block\_size} elements.
For \figurename~\ref{pic:nbody-iters}, \figurename~\ref{pic:nbody-eps}, \figurename~\ref{pic:disolv-iters} and \figurename~\ref{pic:disolv-eps} each SVD reduction of Algorithm \ref{alg:total}
is computed with the relative tolerance $\frac{\tau}{level\_count}$.
We used it instead of $\tau$ to align relative error, reached in
solvation problem on 1 iteration of our algorithm
(see Figure \ref{pic:disolv-iters}).
All the tests, except tests on number of iterations, were performed with 1 iteration of MCBH.
In addition, we tested the optimization of the ranks by the SVD
recompression of the $\mathcal{H}^2$-matrices
\cite{borm-h2book-2010} with different tolerances.
As a starting point for the SVD recompression we used the MCBH
approximation with the accuracy
parameter $\tau = \frac{10^{-8}}{level\_count}$ and 1 additional iteration.
Memory requirement in the figures below is the memory to store both \emph{transfer matrices} and
\emph{interaction matrices}.
The error in figures is relative error of only far field approximation in the spectral norm.
It is computed with the help of Propack package [\url{http://sun.stanford.edu/~rmunk/PROPACK/}].

Tests were performed on a server with 2 Intel(R) Xeon(R) CPU E5504 @
2.00GHz processors with 72GB of RAM.
However, only 2 threads were used (this is default number of threads
for MKL).
Python, Cython and MKL are from the Enthought Python Distribution (EPD 7.3-1 ,64-bit) [https://www.enthought.com/].
Python version is ``2.7.3'', Cython version is ``0.16'', MKL version
is ``10.3-1''.

\subsection{Experiment with particles and interaction matrix \eqref{eq:nbodyaij}}
\label{sec-6-2}

\begin{figure}[H]
\centering
\resizebox{14cm}{!}{\input{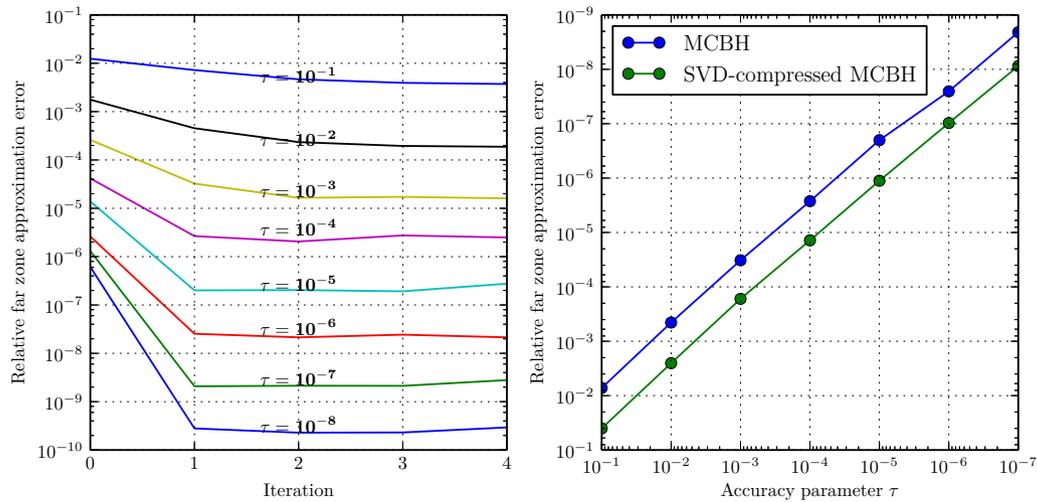}}
\caption{Dependence of approximation error on number of iterations
  (left) and accuracy parameter $\tau$ (right) for the electrostatics
  problem \eqref{eq:nbodyaij}, $N = 100000$.}
\label{pic:nbody-iters}
\end{figure}

\begin{figure}[H]
\centering
\resizebox{14cm}{!}{\input{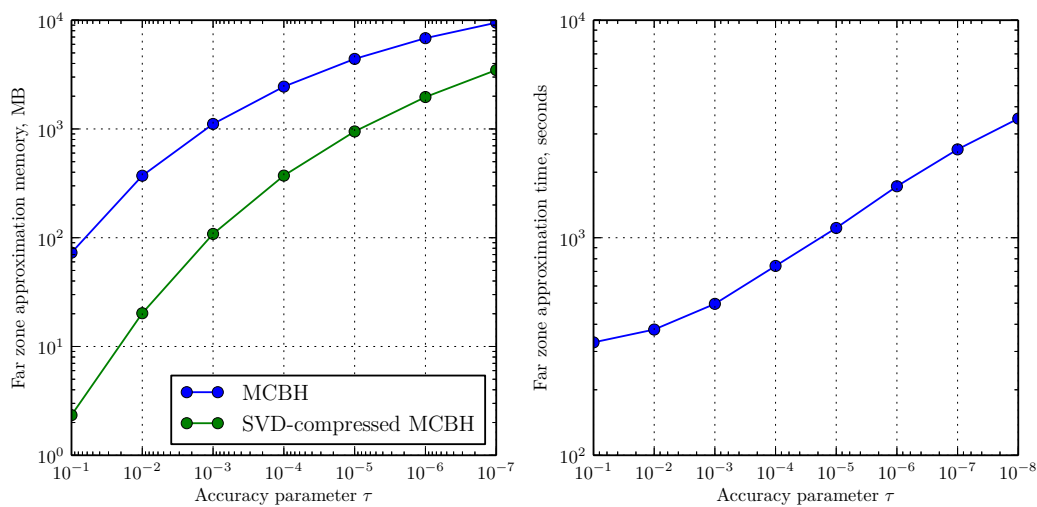}}
\caption{Dependence of approximation time (right) and
  memory (left) on accuracy parameter $\tau$ for the electrostatics
  problem, $N = 100000$.}
\label{pic:nbody-eps}
\end{figure}

Figure \ref{pic:nbody-iters} shows dependence of relative error
on number of iterations and accuracy parameter.
Figure \ref{pic:nbody-eps} shows dependence of approximation time
and memory on accuracy parameter.
To show dependence of approximation time and memory on number of
particles and accuracy parameter, we measured maximum and mean values
for 100 different random particle configurations.
Figures \ref{pic:nbody-max-time}, \ref{pic:nbody-mean-time}, \ref{pic:nbody-max-mem}, \ref{pic:nbody-mean-mem}, \ref{pic:nbody-worst-ratio}, \ref{pic:nbody-mean-ratio} tests used $\tau$ accuracy parameter instead of mentioned $\frac{\tau}{level\_count}$.
Figure \ref{pic:nbody-max-time} and
Figure \ref{pic:nbody-mean-time} show dependence of maximum and
mean approximation time on accuracy parameter and number of particles.
Figure \ref{pic:nbody-max-mem} and
Figure \ref{pic:nbody-mean-mem} show dependence of maximum and
mean approximation memory on accuracy parameter and number of
particles.
Figures \ref{pic:nbody-worst-ratio} and \ref{pic:nbody-mean-ratio} show compress ratio (approximation size vs matrix size).
\begin{figure}[H]
\centering
\resizebox{14cm}{!}{\input{max_factorize_time.pgf}}
\caption{Dependence of maximum (100 tests) approximation time on accuracy parameter (left)
  and number of particles (right) for the electrostatics problem \eqref{eq:nbodyaij}.}
\label{pic:nbody-max-time}
\end{figure}
\begin{figure}[H]
\centering
\resizebox{14cm}{!}{\input{mean_factorize_time.pgf}}
\caption{Dependence of mean (100 tests) approximation time on accuracy parameter (left)
  and number of particles (right) for the electrostatics problem \eqref{eq:nbodyaij}.}
\label{pic:nbody-mean-time}
\end{figure}
\begin{figure}[H]
\centering
\resizebox{14cm}{!}{\input{max_factor_nbytes.pgf}}
\caption{Dependence of maximum (100 tests) approximation memory on accuracy parameter (left)
  and number of particles (right) for the electrostatics problem \eqref{eq:nbodyaij}.}
\label{pic:nbody-max-mem}
\end{figure}
\begin{figure}[H]
\centering
\resizebox{14cm}{!}{\input{mean_factor_nbytes.pgf}}
\caption{Dependence of mean (100 tests) approximation memory on accuracy parameter (left)
  and number of particles (right) for the electrostatics problem \eqref{eq:nbodyaij}.}
\label{pic:nbody-mean-mem}
\end{figure}
\begin{figure}[H]
\centering
\resizebox{14cm}{!}{\input{worst_compress_ratio.pgf}}
\caption{Dependence of worst (100 tests) compress ratio on accuracy parameter (left) and number of particles (right) for the electrostatics problem \eqref{eq:nbodyaij}.}
\label{pic:nbody-worst-ratio}
\end{figure}
\begin{figure}[H]
\centering
\resizebox{14cm}{!}{\input{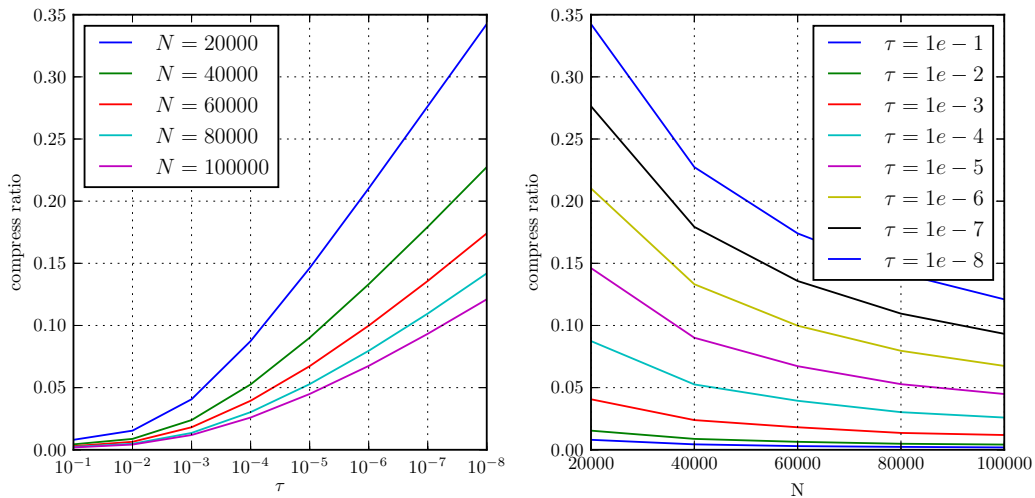}}
\caption{Dependence of mean (100 tests) compress ratio on accuracy parameter (left) and number of particles (right) for the electrostatics problem \eqref{eq:nbodyaij}.}
\label{pic:nbody-mean-ratio}
\end{figure}

\subsection{Boundary integral equation experiment}
\label{sec-6-3}
\begin{figure}[H]
\centering
\resizebox{14cm}{!}{\input{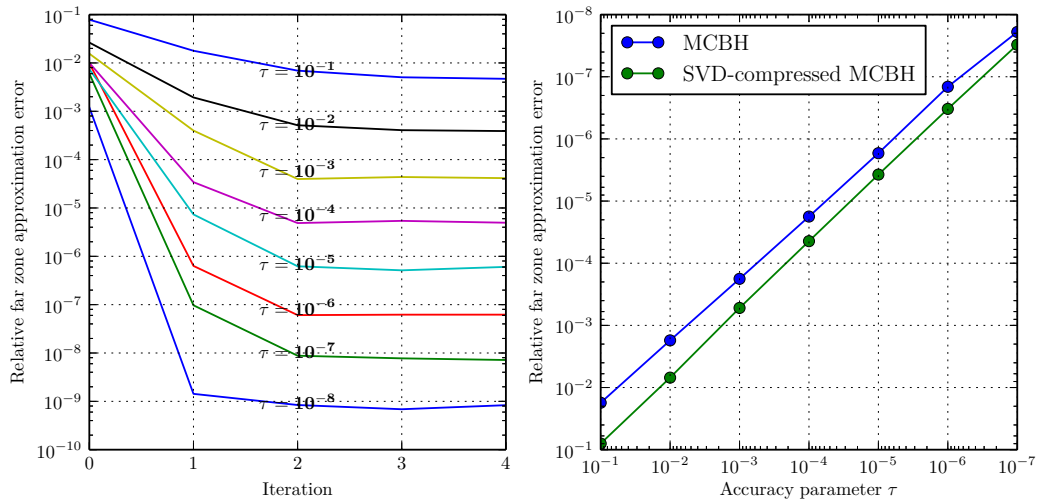}}
\caption{Dependence of approximation error on number of iterations (left) and accuracy parameter $\tau$
(right) for the boundary integral equation \eqref{eq:sigmar}, surface consists
of 222762 discrete elements.}
\label{pic:disolv-iters}
\end{figure}

\begin{figure}[H]
\centering
\resizebox{14cm}{!}{\input{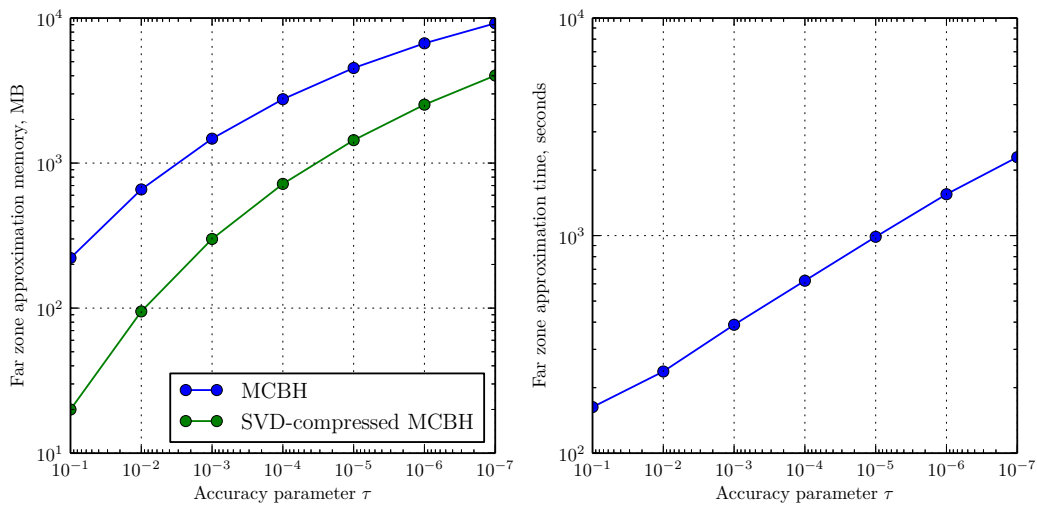}}
\caption{Dependence of approximation time (right) and
  memory (left) on accuracy parameter $\tau$ for the boundary
  integral operator \eqref{eq:sigmar}, surface consists of 222762 discrete elements.}
\label{pic:disolv-eps}
\end{figure}

In this experiment the test surface consists of 222762 triangles.
Figure \ref{pic:disolv-iters} shows dependence of relative error
on number of iterations and accuracy parameter.
Figure \ref{pic:disolv-eps} shows dependence of approximation time
and memory on accuracy parameter.
\section{Comparison with other methods}
\label{sec-7}
    For the comparison we have chosen the methods describe in the paper \cite{bebe-nested-2012}, namely 
    ACAMERGE and ACAGEO. 

We have found that  ACAMERGE method is equivalent to the zero iteration of
proposed multicharge method (MCBH).
As it can be seen from the numerical experiments, this iteration can have limited accuracy. 
Impact of the number of iterations on relative accuracy for different problems is shown in \figurename~\ref{pic:nbody-iters} and \figurename~\ref{pic:disolv-iters}.
So, ACAMERGE cannot get relative error of far field approximation lower
than $10^{-3}$ at least for the double-layer-type problem considered above. 

ACAGEO method is based on the Chebyshev grid introduced in the bounding box of a given cluster node. 
Since the original code from \cite{bebe-nested-2012} is not available, for a fair comparison we have reimplemented
the ACAGEO method. 

\tablename~\ref{table:nbody} and \tablename~\ref{table:pcm} correspond to application of MCBH and ACAGEO to approximation of matrices for two problems from Section~\ref{sec-6}.
Notation used is the following:\\
$\tau$ -- accuracy parameter,\\
$block\_size$ -- maximum size of leaf node (number of elements in
it),\\
$iters$ -- number of iterations for MCBH method,\\
$time$ -- factorization time in seconds,\\
$error$ -- relative spectral error of approximation, measured only by far field,\\
$points$ -- number of Chebyshev points for ACAGEO method (number of
potential basis elements for each node of cluster trees).

\begin{table}[H]
\centering
\caption{Comparison of time and error of MCBH and ACAGEO for electrostatics problem \eqref{eq:nbodyaij}.}
\label{table:nbody}
\begin{tabular}{|c|c|c|c|c|c|c|c|}
\hline
& & \multicolumn{3}{|c|}{MCBH} & \multicolumn{3}{|c|}{ACAGEO} \\ \hline
$\tau$ & block\_size & iters & time & error & points & time & error \\ \hline
\multirow{2}{*}{$10^{-2}$} & \multirow{2}{*}{25} & \multirow{2}{*}{0} & \multirow{2}{*}{9.7} & \multirow{2}{*}{$2.2*10^{-2}$} & 8 & 10.7 & $1.8*10^{-2}$ \\
& & & & & 27 & 29.5 & $2.2*10^{-2}$ \\ \hline
$10^{-3}$ & 25 & 0 & 15.5 & $3.2*10^{-3}$ & 27 & 30.5 & $3.8*10^{-3}$ \\ \hline
$10^{-4}$ & 25 & 0 & 26.8 & $4.3*10^{-4}$ & 27 & 32.3 & $6.3*10^{-4}$ \\ \hline 
\multirow{2}{*}{$10^{-5}$} & \multirow{2}{*}{25} & 0 & 46.3 & $1.2*10^{-4}$ & 27 & 35 & $2.4*10^{-4}$ \\
& & 1 & 148 & $3*10^{-5}$ & 64 & 217.4 & $1.8*10^{-4}$\\ \hline
\multirow{2}{*}{$10^{-6}$} & \multirow{2}{*}{25} & \multirow{2}{*}{1} & \multirow{2}{*}{262} & \multirow{2}{*}{$3.3*10^{-6}$} & 27 & 33.4 & $1.8*10^{-4}$ \\
& & & & & 64 & 221 & $4*10^{-5}$ \\ \hline
\end{tabular}
\end{table}

\begin{table}[H]
\centering
\caption{Comparison of time and error of MCBH and ACAGEO for solvation problem.}
\label{table:pcm}
\begin{tabular}{|c|c|c|c|c|c|c|c|}
\hline
& & \multicolumn{3}{|c|}{MCBH} & \multicolumn{3}{|c|}{ACAGEO} \\ \hline
$\tau$ & block\_size & iters & time & error & points & time & error \\ \hline
\multirow{2}{*}{$10^{-2}$} & 25 & 0 & 11.5 & $8.4*10^{-2}$ & 8 & 13.2 & $2.6*10^{-1}$ \\
& 25 & 1 & 29 & $6.4*10^{-2}$ & 27 & 35.7 & $9.6*10^{-2}$ \\ \hline
\multirow{2}{*}{$10^{-3}$} & 25 & 0 & 18.6 & $2.4*10^{-2}$ & 8 & 13.3 & $2.6*10^{-1}$ \\
& 25 & 1 & 50 & $1.1*10^{-2}$ & 27 & 37 & $3.3*10^{-2}$ \\ \hline
\multirow{3}{*}{$10^{-4}$} & 25 & 0 & 30.8 & $1.2*10^{-2}$ & 8 & 13.2 & $2.6*10^{-1}$ \\
& 25 & 1 & 101.9 & $1.6*10^{-3}$ & 27 & 38.2 & $1.7*10^{-2}$ \\
& 25 & 2 & 174.9 & $1.6*10^{-3}$ & 64 & 163 & $1.5*10^{-2}$ \\ \hline
\multirow{3}{*}{$10^{-5}$} & 25 & 0 & 45.9 & $9*10^{-3}$ & 27 & 38.6 & $1.4*10^{-2}$ \\
& 25 & 1 & 181.6 & $2.1*10^{-4}$ & 64 & 166.3 & $7.6*10^{-3}$ \\
& 25 & 2 & 321.8 & $2*10^{-4}$ & 125 & 609 & $8.6*10^{-3}$\\ \hline
\end{tabular}
\end{table}

\section{Conclusions and future work}
\label{sec-8}
The proposed algorithm is robust in the sense that high approximation accuracies are
achievable using at most 2 additional iterations.  It also does not require the construction of apriori
representing sets. An important question that needs to be solved is the
existence of a good MCBH-type approximation, provided that a good
$\mathcal{H}^2$-approximation exists. A preliminary study shows that it
is possible to derive the existence of a quasi-optimal approximation
in the spirit of the recent paper \cite{sav-qott-2013pre}. 

Computational speed and memory requirements of the algorithm can be
improved in several ways. The algorithm can be readily parallelized
since it has the usual FMM structure. Also, in the case when the
elements of a matrix can be computed in a cheap way, the interaction matrices can be computed
online, and the memory consumption becomes much smaller. For the case
when the computation of a single element of a matrix is expensive (for example, in
the Galerkin method for the solution of a boundary integral equation), the application
of the $\mathcal{H}^2$-recompression algorithm may give a noticeable
reduction in the storage. 

The code is available as a part of the h2tools package 
\url{http://bitbucket.org/muxas/h2tools}, and we plan to extended it with
different compression algorithms and solvers. It will be interesting to apply the new method to
different problems, for example to the problem of vortex ring dynamics
 \cite{stavtsev-crossvortex-2012}.

\bibliography{mcbh-new,bibtex/our}
\bibliographystyle{wileyj.bst}
\end{document}

%% file: include.tex
\tikzset{full/.style = {fill = gray}, lowrank/.style = {}, fullrank/.style = {fill = gray}}
\newsavebox{\hmatrixA}
\savebox{\hmatrixA}{%
\begin{tikzpicture}[y = -1cm, overlay]
\draw[full] (0, 0) rectangle +(16, 16);
\end{tikzpicture}%
}

\newsavebox{\hmatrixB}
\savebox{\hmatrixB}{%
\begin{tikzpicture}[y = -1cm, overlay, remember picture]
\draw (0, 0) rectangle +(16, 16);
\draw[full] (0, 0) rectangle +(8,8);
\draw[full] (8, 8) rectangle +(8,8);
\end{tikzpicture}%
}

\newsavebox{\hmatrixC}
\savebox{\hmatrixC}{%
\begin{tikzpicture}[y = -1cm, overlay, remember picture]
\draw (0, 0) rectangle +(16, 16);
\draw (0, 0) rectangle +(8,8);
\draw (8, 8) rectangle +(8,8);
\foreach \i in {0, 1, 2, 3} \filldraw[full] (4*\i, 4*\i) rectangle +(4, 4);
\end{tikzpicture}%
}

\newsavebox{\hmatrixD}
\savebox{\hmatrixD}{%
\begin{tikzpicture}[y = -1cm, overlay, remember picture]
\draw (0, 0) rectangle +(16, 16);
\draw (0, 0) rectangle +(8, 8);
\draw (8, 8) rectangle +(8, 8);
\foreach \i in {0, 1, 2, 3} \draw (4*\i, 4*\i) rectangle +(4, 4);
\foreach \i in {0, 1, ..., 7} \filldraw[full] (2*\i, 2*\i) rectangle +(2, 2);
\end{tikzpicture}%
}

\newsavebox{\hmatrixE}
\savebox{\hmatrixE}{%
\begin{tikzpicture}[y = -1cm, overlay, remember picture]
\draw (0, 0) rectangle +(16, 16);
\draw (0, 0) rectangle +(8, 8);
\draw (8, 8) rectangle +(8, 8);
\foreach \i in {0, 1, 2, 3} \draw (4*\i, 4*\i) rectangle +(4, 4);
\foreach \i in {0, 1, ..., 7} \draw (2*\i, 2*\i) rectangle +(2, 2);
\foreach \i in {0, 1, ..., 15} \filldraw[full] (\i, \i) rectangle +(1, 1);
\end{tikzpicture}%
}


\newsavebox{\hmatrixF}
\savebox{\hmatrixF}{%
\begin{tikzpicture}[y = -1cm, overlay, remember picture]
\filldraw[full] (0, 0) rectangle +(8, 8);
\filldraw[full] (0, 8) rectangle +(8, 8);
\filldraw[full] (8, 0) rectangle +(8, 8);
\filldraw[full] (8, 8) rectangle +(8, 8);
\end{tikzpicture}%
}

\newsavebox{\hmatrixG}
\savebox{\hmatrixG}{%
\begin{tikzpicture}[y = -1cm, overlay, remember picture]
\draw (0, 0) rectangle +(16, 16);
\foreach \i in {0, 1, 2, 3} \foreach \j in {0, 1, 2, 3} \filldraw[full] (4*\i, 4*\j) rectangle +(4, 4);
\foreach \i in {0, 1, 2, 3} \draw[fill = white] (12-4*\i, 4*\i) rectangle +(4, 4);
\end{tikzpicture}%
}

\newsavebox{\hmatrixH}
\savebox{\hmatrixH}{%
\begin{tikzpicture}[y = -1cm, overlay, remember picture]
\foreach \i in {(12, 0), (8, 4), (4, 8), (0, 12)} \draw[fill = white] \i rectangle +(4, 4);
\foreach \i in {(6, 0), (8, 0), (10, 0), (4, 2), (10, 2), (2, 4), (12, 4), (14, 4), (0, 6), (14, 6), (0, 8), (14, 8), (0, 10), (2, 10), (12, 10), (4, 12), (10, 12), (4, 14), (6, 14), (8, 14)} \draw[fill = white] \i rectangle +(2, 2);
\foreach \i in {(0, 0), (2, 0), (4, 0), (0, 2), (2, 2), (6, 2), (8, 2), (0, 4), (4, 4), (6, 4), (2, 6), (4, 6), (6, 6), (12, 6), (2, 8), (8, 8), (10, 8), (12, 8), (8, 10), (10, 10), (14, 10), (6, 12), (8, 12), (12, 12), (14, 12), (10, 14), (12, 14), (14, 14)} \draw[full] \i rectangle +(2, 2);
\end{tikzpicture}%
}

\newsavebox{\hmatrixI}
\savebox{\hmatrixI}{%
\begin{tikzpicture}[y = -1cm, overlay, remember picture]
\foreach \i in {(12, 0), (8, 4), (4, 8), (0, 12)} \draw[fill = white] \i rectangle +(4, 4);
\foreach \i in {(6, 0), (8, 0), (10, 0), (4, 2), (10, 2), (2, 4), (12, 4), (14, 4), (0, 6), (14, 6), (0, 8), (14, 8), (0, 10), (2, 10), (12, 10), (4, 12), (10, 12), (4, 14), (6, 14), (8, 14)} \draw[fill = white] \i rectangle +(2, 2);
\foreach \i in {(3, 0), (4, 0), (5, 0), (2, 1), (5, 1), (1, 2), (6, 2), (7, 2), (9, 2), (0, 3), (7, 3), (8, 3), (0, 4), (7, 4), (0, 5), (1, 5), (6, 5), (2, 6), (5, 6), (13, 6), (2, 7), (3, 7), (4, 7), (12, 7), (3, 8), (11, 8), (12, 8), (13, 8), (2, 9), (10, 9), (13, 9), (9, 10), (14, 10), (15, 10), (8, 11), (15, 11), (7, 12), (8, 12), (15, 12), (6, 13), (8, 13), (9, 13), (14, 13), (10, 14), (13, 14), (10, 15), (11, 15), (12, 15)} \draw[fill = white] \i rectangle +(1, 1);
\foreach \i in {(0, 0), (1, 0), (2, 0), (0, 1), (1, 1), (3, 1), (4, 1), (0, 2), (2, 2), (3, 2), (8, 2), (1, 3), (2, 3), (3, 3), (6, 3), (9, 3), (1, 4), (4, 4), (5, 4), (6, 4), (4, 5), (5, 5), (7, 5), (3, 6), (4, 6), (6, 6), (7, 6), (12, 6), (5, 7), (6, 7), (7, 7), (13, 7), (2, 8), (8, 8), (9, 8), (10, 8), (3, 9), (8, 9), (9, 9), (11, 9), (12, 9), (8, 10), (10, 10), (11, 10), (9, 11), (10, 11), (11, 11), (14, 11), (6, 12), (9, 12), (12, 12), (13, 12), (14, 12), (7, 13), (12, 13), (13, 13), (15, 13), (11, 14), (12, 14), (14, 14), (15, 14), (13, 15), (14, 15), (15, 15)} \draw[full] \i rectangle +(1, 1);
\end{tikzpicture}%
}

%% file: hmatrix.tex
\begin{tikzpicture}[y = -1cm, inner sep = 0]
\begin{scope}[every node/.style = {scale = 1/8}]
\path (0, 0) rectangle +(2, 2) node at +(0, 0) {\usebox{\hmatrixA}};
\path (3, 0) rectangle +(2, 2) node at +(0, 0) {\usebox{\hmatrixB}};
\path (6, 0) rectangle +(2, 2) node at +(0, 0) {\usebox{\hmatrixC}};
\path (9, 0) rectangle +(2, 2) node at +(0, 0) {\usebox{\hmatrixD}};
\path (12, 0) rectangle +(2, 2) node at +(0, 0) {\usebox{\hmatrixE}};
\end{scope}
\draw [->] (2+0.2, 1) -- +(0.6, 0);
\draw [->] (5+0.2, 1) -- +(0.6, 0);
\draw [->] (8+0.2, 1) -- +(0.6, 0);
\draw [->] (11+0.2, 1) -- +(0.6, 0);
\begin{scope}[shift = {(0, 3)}]
\begin{scope}[every node/.style = {scale = 1/8}]
\path (0, 0) rectangle +(2, 2) node at +(0, 0) {\usebox{\hmatrixA}};
\path (3, 0) rectangle +(2, 2) node at +(0, 0) {\usebox{\hmatrixF}};
\path (6, 0) rectangle +(2, 2) node at +(0, 0) {\usebox{\hmatrixG}};
\path (9, 0) rectangle +(2, 2) node at +(0, 0) {\usebox{\hmatrixH}};
\path (12, 0) rectangle +(2, 2) node at +(0, 0) {\usebox{\hmatrixI}};
\end{scope}
\draw [->] (2+0.2, 1) -- +(0.6, 0);
\draw [->] (5+0.2, 1) -- +(0.6, 0);
\draw [->] (8+0.2, 1) -- +(0.6, 0);
\draw [->] (11+0.2, 1) -- +(0.6, 0);
\end{scope}
\begin{scope}[shift = {(16, 1)}]
\draw (0, 0) rectangle +(1, 1) node [midway, xshift = 4] {\footnotesize admissable far block};
\draw[full] (0, 2) rectangle +(1, 1) node [midway, xshift = 4]
{\footnotesize admissable close block};
\end{scope}
\end{tikzpicture}

%% file: blockstring.tex
\begin{tikzpicture}[y = -1cm, inner sep = 0, scale = 1/4, every node/.style = {scale = 1/4}]
\path (0, 0) rectangle +(16, 16) node at +(0, 0) {\usebox{\hmatrixI}};
\draw [very thin, decorate,decoration={brace, mirror}, xshift = -4] (0, 4) -- +(0, 2) node [midway,xshift=-20, scale = 4] {\footnotesize $t$};
\draw [decorate, decoration = {brace}, yshift = 4] (2, 0) -- +(2, 0) node [midway, yshift = 20, scale = 4] {\footnotesize $s_0$};
\draw [decorate, decoration = {brace}, yshift = 4] (12, 0) -- +(2, 0) node [midway, yshift = 20, scale = 4] {\footnotesize $s_1$};
\draw [decorate, decoration = {brace}, yshift = 4] (14, 0) -- +(2, 0) node [midway, yshift = 20, scale = 4] {\footnotesize $s_2$};
\draw[fill = gray!50] (2, 4) rectangle +(2, 2);
\draw[fill = gray!50] (12, 4) rectangle +(2, 2);
\draw[fill = gray!50] (14, 4) rectangle +(2, 2);
\draw[fill = gray!50] (0, 18) rectangle +(2, 2);
\draw[pattern = north west lines] (8, 4) rectangle +(4, 2);
\draw[pattern = north west lines] (8, 18) rectangle +(2, 2);
\node[scale = 4, anchor = north west] at (2.5, 18.5) {\footnotesize $R^*_t$};
\node[scale = 4, anchor = north west] at (10.5, 18.5) {\footnotesize $R^+_t$};

\begin{scope}[shift = {(30, 0)}]
\path (0, 0) rectangle +(16, 16) node at +(0, 0) {\usebox{\hmatrixI}};
\draw [very thin, decorate,decoration={brace, mirror}, xshift = -4] (0, 4) -- +(0, 2) node [midway,xshift=-20, scale = 4] {\footnotesize $t$};
\draw[fill = gray!50] (2, 4) rectangle +(2, 2);
\draw[fill = gray!50] (12, 4) rectangle +(2, 2);
\draw[fill = gray!50] (14, 4) rectangle +(2, 2);
\draw[fill = gray!50] (8, 18) rectangle +(2, 2);
\draw[fill = gray!50] (8, 4) rectangle +(4, 2);
\node[scale = 4, anchor = north west] at (10.5, 18.5) {\footnotesize $R_t$};
\end{scope}
\end{tikzpicture}%